\documentclass{article}%
\usepackage{amssymb}
\usepackage{amsfonts}
\usepackage{graphicx}
\usepackage{amsmath}%
\providecommand{\U}[1]{\protect\rule{.1in}{.1in}}
\newtheorem{theorem}{Theorem}

\newtheorem{corollary}[theorem]{Corollary}

\newtheorem{lemma}[theorem]{Lemma}

\newtheorem{proposition}[theorem]{Proposition}
\newtheorem{remark}[theorem]{Remark}

\begin{document}

\begin{center}
{\LARGE Central limit theorem started at a point for additive functionals of
reversible Markov chains}
\end{center}

\bigskip

\textbf{Abbreviated Title: CLT for reversible Markov chains}

\bigskip

\begin{center}
Christophe Cuny and Magda Peligrad\footnote{Supported in part by a Charles
Phelps Taft Memorial Fund grant and NSA\ grant H98230-09-1-0005
\par
{}}
\end{center}

\bigskip

Equipe ERIM, University of New Caledonia, B.P. R4, 98800 Noum\'{e}a, New
Caledonia. E-mail address: cuny@univ-nc.nc \ 

Department of Mathematical Sciences, University of Cincinnati, PO Box 210025,
Cincinnati, Oh 45221-0025, USA. E-mail address: peligrm@ucmail.uc.edu

\begin{center}
\bigskip

Abstract
\end{center}

In this paper we study the almost sure central limit theorem started from a
point for additive functionals of a stationary and ergodic Markov chain via a
martingale approximation in the almost sure sense. As a consequence we derive
the quenched CLT\ for additive functionals or reversible Markov chains.

\section{Introduction}

\ \ \ \ In 1986 Kipnis and Varadhan proved the functional form of the central
limit theorem for additive functionals of stationary reversible ergodic Markov
chains under a minimal spectral assumption. This result was established with
respect to the stationary probability law of the chain. In their remark (1.7)
Kipnis and Varadhan raised the question if their result also holds with
respect to the law of the Markov chain started from $x$, for almost all $x$.

The central limit theorem for stationary Markov chains with normal operator
holds under a similar spectral assumption as discovered by Gordin and Lifshitz
(1981) (see also and Borodin and Ibragimov, 1994, ch 4 sections 7-8,
Derriennic and Lin, 1996, Zhao and Woodroofe, 2008).

The problem of quenched CLT for normal stationary and ergodic Markov chains
was considered by Derriennic and Lin (2001 a and b) and Cuny (2009 a) under
some reinforced assumptions. The results concerning reversible Markov chains
are usually obtained as corollaries. Our Theorem \ref{REV}, provides a central
limit theorem started at a point for reversible Markov chains that improves
the results known in the literature and answers a question raised by
Derriennic and Lin (2001 a, section 7). It is specific for reversible Markov
chains and the result does not hold for normal Markov chains without further
reinforcing the spectral condition (see Proposition \ref{P1} below).

Examples of reversible Markov chains frequently appear in the study of
infinite systems of particles, random walks or processes in random media. In
this context, the almost sure validity of the central limit theorem refers to
as "quenched" media as opposed to the "annealed" media. A simple example of
Normal Markov chain is a random walk on a compact group.

Our method of proof follows the classical line of martingale approximation,
stressing now the study of the almost sure asymptotic results for the rest.

The paper is organized as following. In Section 2 we give the definitions,
recall some needed results and formulate our result. Section 3 contains the
proof of the theorem and also general results on the quenched CLT for
stationary processes that have interest in themselves. Applications to random
walks on compact groups are given in Section 4. A general maximal inequality
is given in the Appendix.

\section{Definitions Background and Results}

\subsection{\textbf{Notations}}

\ \ \ \ We assume that $(\xi_{n})_{n\in\mathbb{Z}}$ is a stationary and
ergodic Markov chain defined on a probability space $(\Omega,\mathcal{F}%
,\mathbb{P})$ with values in a measurable space $(S,\mathcal{A})$. The
marginal distribution is denoted by $\pi(A)=P(\xi_{0}\in A)$ and we suppose
there is a regular conditional distribution for $\xi_{1}$ given $\xi_{0}$
denoted by $Q(x,A)=\mathbb{P}(\xi_{1}\in A|\,\xi_{0}=x)$. In addition $Q$
denotes the operator {acting via $(Qf)(x)=\int_{S}f(s)Q(x,ds).$ Next let
$\mathbb{L}_{0}^{2}(\pi)$ be the set of functions on $S$ such that $\int
f^{2}d\pi<\infty$ and $\int fd\pi=0.$ For a function }${f}\in${$\mathbb{L}%
_{0}^{2}(\pi)$ let $X_{i}=f(\xi_{i})$, $S_{n}=\sum\limits_{i=1}^{n}X_{i}$.
Denote by $\mathcal{F}_{k}$ the $\sigma$--field generated by $\xi_{i}$ with
$i\leq k.$ For any integrable variable $X$ we denote $\mathbb{E}%
_{k}(X)=\mathbb{E}(X|\mathcal{F}_{k}).$ In our notation $\mathbb{E}_{0}%
(X_{1})=(Qf)(\xi_{0})=\mathbb{E}(X_{1}|\xi_{0}).$ We denote by }${{||X||}_{2}%
}$ the norm in {$\mathbb{L}_{2}$}$(\Omega,\mathcal{F},\mathbb{P})$ and by
{$||f||_{2}$ the norm in $\mathbb{L}_{2}(\pi).$}

Denote by $\mathbb{P}^{x}$ and $\mathbb{E}^{x}$ the regular probability and
conditional expectation given $X_{0}=x$. By the CLT started at a point
(quenched) we understand that the CLT holds for $\pi-$almost all $x\in S,$
under the measure $\mathbb{P}^{x}$.

The Markov chain is called normal if $QQ^{\ast}=Q^{\ast}Q$ on {$\mathbb{L}%
_{2}(\pi)$, }where the adjoint operator $Q^{\ast}$ is defined by
$<Qf,g>=<f,Q^{\ast}g>$,$\ $for{\ every }${f}${\ and }${g}${\ in $\mathbb{L}%
_{2}(\pi)$}$.$ {For every }${f}${\ in $\mathbb{L}^{2}(\pi)$ we denote by
}$\rho_{f}$ the spectral measure of $f$ with respect to $Q$ on the unit disk
$D.$

The Markov chain is called reversible if, $Q=Q^{\ast}.$ For reversible Markov
chains the spectral measure is concentrated on $[-1,1]$.

\subsection{ \textbf{Background and Results}}

\ \ \ \ The central limit theorem for normal Markov chains was announced by
Gordin and Lifshitz, (1981) and also treated by Derriennic and Lin (1996). Its
proof can be found in Borodin and Ibragimov (1994, ch. 4, sections 7). The
result basically states that for every stationary ergodic Markov chain with
$Q$ normal and $f\in${$\mathbb{L}_{0}^{2}(\pi)$ satisfying}%
\begin{equation}
\int\nolimits_{D}\frac{1}{|1-z|}\rho_{f}(dz)<\infty\text{ ,} \label{SN}%
\end{equation}
then the central limit theorem holds for $S_{n}/\sqrt{n}.$

For reversible Markov chains, under the same spectral condition, Kipnis and
Varadhan (1986) proved that the central limit theorem holds in its functional
form. Let us introduce the notation $[x]$ as being the integer part of $x$ and
for $0\leq t\leq1$ define%
\[
W_{n}(t)=\frac{S_{[nt]}}{\sqrt{n}}\text{ .}%
\]
$W_{n}(t)$ belongs to the space $D[0,1]$ of functions continuous from the
right with limits from the left endowed with the uniform topology.

\begin{theorem}
\label{B}(Kipnis and Varadhan, 1986). Assume $(\xi_{n})_{n\in\mathbb{Z}}$ is a
stationary ergodic reversible Markov chain and $f\in${$\mathbb{L}_{0}^{2}%
(\pi)$ satisfies}%
\begin{equation}
\int\nolimits_{-1}^{1}\frac{1}{1-t}\rho_{f}(dt)<\infty\text{ .} \label{SR}%
\end{equation}
Then, $W_{n}(t)$ converges weakly to $|\sigma_{f}|W$, where%
\[
\sigma_{f}^{2}=\int\nolimits_{-1}^{1}\frac{1+t}{1-t}\rho_{f}(dt)\text{ }%
\]
and $W(t)$ denotes the standard Brownian motion on $[0,1].$
\end{theorem}

Spectral condition (\ref{SR}) (see Kipnis and Varadhan, 1986) is equivalent
to
\begin{equation}
\lim_{n\rightarrow\infty}\frac{var(S_{n})}{n}=\sigma_{f}^{2}\text{ .}
\label{var}%
\end{equation}

An important problem is to investigate the validity of the almost sure central
limit theorem started at a point for additive functionals of stationary
ergodic normal Markov chains or reversible Markov chains under condition
(\ref{SN}). For normal Markov chains Derriennic and Lin (2001, p.515) pointed
out that in general condition (\ref{SN}) does not assure the validity of the
CLT started at a point. In the context of normal stationary and ergodic Markov
chains the problem of quenched CLT was considered by Derriennic and Lin (1996,
2001 a) and Cuny (2009) under some reinforced assumptions.

For reversible Markov chains, we shall improve their results.

\begin{theorem}
\label{REV}Assume $(\xi_{i})_{i\in\mathbb{Z}}$ is a stationary and ergodic
reversible Markov chain and $f\in\mathbb{L}_{0}^{2}(\pi)$ satisfies
\begin{equation}
\int\nolimits_{-1}^{1}\frac{(\log^{+}|\log(1-t)|)^{2}}{1-t}\rho_{f}%
(dt)<\infty\text{ .} \label{SR2}%
\end{equation}
Then, there is a martingale with stationary and ergodic differences
$(M_{n})_{n\geq1}$ adapted to $(${$\mathcal{F}_{n})_{n\geq1}$ }such that for
almost all $x\in S$
\begin{equation}
\mathbb{E}^{x}\frac{[S_{n}-M_{n}]^{2}}{n}\rightarrow0\text{ ,} \label{LimR}%
\end{equation}
and for almost all $x\in S$%
\[
\frac{S_{n}}{\sqrt{n}}\Longrightarrow|\sigma_{f}|N\text{ under }\mathbb{P}%
^{x}\text{ .}%
\]
Here $N$ denotes a standard normal variable.
\end{theorem}

As a matter of fact we prove the martingale approximation (\ref{LimR}) in the
following equivalent formulation:

\begin{remark}
Under conditions of Theorem \ref{REV} there is a martingale with stationary
and ergodic differences $(M_{n})_{n\geq1}$ adapted to $(${$\mathcal{F}%
_{n})_{n\geq1}$ }such that
\begin{equation}
\mathbb{E}_{0}\frac{[S_{n}-M_{n}]^{2}}{n}\rightarrow0\text{ almost surely and
in }\mathbb{L}_{1}\text{ .} \label{LIMR}%
\end{equation}

\end{remark}

We point out that condition (\ref{SR2}) of Theorem \ref{REV} does not assure
the validity of the CLT started at a point for normal Markov chains. Combining
the example in Derriennic and Lin (2001 a, p. 515) with Cuny (2009 b,
Proposition 2.4 and the last part of Theorem 2.5, with $b(n)^{2}=1/\log n$) we conclude:

\begin{proposition}
\label{P1}Let $\varepsilon\in[0,1[$. There is a stationary and ergodic normal
Markov chain and a function $f$ such that
\begin{equation}
\int\nolimits_{D}\frac{|\log(|1-z|)|}{|1-z|}\rho_{f}(dz)<\infty\label{SN1}%
\end{equation}
and%
\[
\lim\sup_{n\rightarrow\infty}\frac{S_{n}}{\sqrt{n (\log\log n)^{\varepsilon}}
}=\infty\text{ a.s.}%
\]

\end{proposition}

Everywhere in the paper we denote by%
\begin{equation}
V_{n}(f)=(I+Q+...+Q^{n-1})(f) \label{DefV}%
\end{equation}
and
\begin{equation}
H_{n}(x,y)=V_{n}(f)(y)-QV_{n}(f)(x)\text{ .} \label{Def H}%
\end{equation}
Define the projector operator
\begin{equation}
\mathbb{P}_{j}(X)=\mathbb{E}_{j}(X)-\mathbb{E}_{j-1}(X)\text{ .}
\label{Defproj}%
\end{equation}
Then, in terms of $V_{n},$ and using the Markov property
\[
\mathbb{P}_{1}(S_{n})=V_{n}(\xi_{1})-QV_{n}(\xi_{0})=H_{n}(\xi_{0},\xi
_{1})\text{ .}%
\]

\begin{remark}
It is interesting to notice that by Cuny (2009 a, Lemma 2.1) condition
(\ref{SR2}) is equivalent to%
\[
\sum\nolimits_{n}\frac{(\log\log n)^{2}||V_{n}(f)||_{2}^{2}}{n^{2}}%
<\infty\text{ .}%
\]

\end{remark}

\section{Preparatory results}

\subsection{Quenched CLT\ for stationary and ergodic sequences}

We establish first the central limit theorem started at a point for general
additive functionals of Markov chains. Since any stationary sequence
$(X_{k})_{k\in\mathbb{Z}}$ can be viewed as a function of a Markov process
$\xi_{k}=(X_{i};i\leq k),$ for the function $g(\xi_{k})=X_{k}$ the results in
this section are valid for any stationary and ergodic sequence of random
variables. Below and everywhere in the paper a.s. denotes either convergence
or equality in the almost sure sense.

\begin{proposition}
\label{LimAS}Assume $(\xi_{n})_{n\in\mathbb{Z}}$ is a stationary ergodic
Markov chain and $f\in${$\mathbb{L}_{0}^{2}(\pi)$ satisfies}%
\begin{equation}
\frac{\mathbb{E}_{0}(S_{n})}{\sqrt{n}}\rightarrow0\text{ a.s. ,} \label{ASDL}%
\end{equation}%
\begin{equation}
H_{m}(\xi_{1},\xi_{0})\rightarrow H(\xi_{1},\xi_{0})\text{ converges a.s. and
\ }\mathbb{E[}\sup_{m}H_{m}^{2}(\xi_{1},\xi_{0})]<\infty\text{ .} \label{AS}%
\end{equation}
Then the conclusion of Theorem \ref{REV} holds.
\end{proposition}

\textbf{Proof}.

Starting from condition (\ref{AS}) we notice that this implies $H_{m}(\xi
_{1},\xi_{0})\rightarrow H(\xi_{1},\xi_{0})$ in $\mathbb{L}_{2}.$ Since
$\mathbb{E}_{0}[H_{m}(\xi_{1},\xi_{0})]=0$ a.s. we conclude that
$\mathbb{E}_{0}[H(\xi_{1},\xi_{0})]=0$ a.s. We consider the sequence of
martingale differences $[H(\xi_{k-1},\xi_{k})]_{k\geq1}$. Then, we define a
martingale
\[
M_{n}=\sum_{k=1}^{n}H(\xi_{k-1},\xi_{k})
\]
adapted to $\mathcal{F}_{n}$ and we shall approximate $S_{n}$ by $M_{n}.$ We
shall use the decomposition:%
\[
\frac{S_{n}-M_{n}}{\sqrt{n}}=\frac{1}{\sqrt{n}}\mathbb{E}_{0}(S_{n})+\frac
{1}{\sqrt{n}}[S_{n}-\mathbb{E}_{0}(S_{n})-M_{n}]\text{ .}%
\]
In order to establish (\ref{LimR}), by (\ref{ASDL}), it remains to show that
\[
\frac{1}{n}\mathbb{E}_{0}[S_{n}-\mathbb{E}_{0}(S_{n})-M_{n}]^{2}%
\rightarrow0\text{ a.s.}%
\]
We use now a traditional decomposition of $S_{n}$ in martingale differences by
using the projections on consecutive sigma algebras:%
\[
S_{n}-\mathbb{E}_{0}(S_{n})=[S_{n}-\mathbb{E}_{n-1}(S_{n})]+[\mathbb{E}%
_{n-1}(S_{n})-\mathbb{E}_{n-2}(S_{n})]+...+[\mathbb{E}_{1}(S_{n}%
)-\mathbb{E}_{0}(S_{n})]\text{ .}%
\]
So, by using notation (\ref{Def H}) we have the martingale decomposition
\[
\frac{1}{\sqrt{n}}[S_{n}-\mathbb{E}_{0}(S_{n})-M_{n}]=\frac{1}{\sqrt{n}}%
\sum_{k=1}^{n}[H_{n-k+1}(\xi_{k-1},\xi_{k})-H(\xi_{k-1},\xi_{k})]\text{ .}%
\]
Then, by the martingale properties followed by the change of variable and
using the notation
\begin{equation}
H_{m}-H=G_{m}\text{ ,}\label{defK}%
\end{equation}
it follows that
\begin{gather*}
\frac{1}{n}\mathbb{E}_{0}[S_{n}-\mathbb{E}_{0}(S_{n})-M_{n}]^{2}=\frac{1}%
{n}\sum_{k=1}^{n}\mathbb{E}_{0}[G_{n-k+1}(\xi_{k-1},\xi_{k})]^{2}\\
=\frac{1}{n}\sum_{j=1}^{n}\mathbb{E}_{0}[G_{j}(\xi_{n-j},\xi_{n-j+1}%
)]^{2}\text{ .}%
\end{gather*}
Let $N$ be fixed. For $n$ sufficiently large we decompose the last sum into a
sum from $1$ to $N$ and one from $N+1$ to $n.$ Then,%
\begin{gather*}
\frac{1}{n}\mathbb{E}_{0}[S_{n}-\mathbb{E}_{0}(S_{n})-M_{n}]^{2}\leq\frac
{1}{n}\sum_{j=1}^{N}\mathbb{E}_{0}[G_{j}(\xi_{n-j},\xi_{n-j+1})]^{2}+\\
\frac{1}{n}\sum_{j=N+1}^{n}\mathbb{E}_{0}[G_{j}(\xi_{n-j},\xi_{n-j+1}%
)]^{2}=I_{n}(N)+II_{n}(N)\text{ .}%
\end{gather*}
We shall prove the almost sure negligibility of each term in the last sum.

\smallskip

We treat first $I_{n}(N).$ Define $L_{j}(x):=\int_{S}G_{j}^{2}(x,y)Q(x,dy)$.
Then
\[
\mathbb{E}_{0}[G_{j}(\xi_{n-j},\xi_{n-j+1})]^{2}=Q^{n-j}L_{j}(\xi_{0})\text{
.}%
\]
Now, by the ergodic theorem for the Markov operator $Q$,
\[
\frac{Q^{n-j}L_{j}(\xi_{0})}{n}=\frac{n}{n-j}\frac{Q^{n-j}L_{j}(\xi_{0})}%
{n-j}\rightarrow0~\text{a.s.}%
\]
for every $1\leq j\leq N$, and $I_{n}(N)\rightarrow0~\text{a.s.}$.



\smallskip

Now we treat $II_{n}(N)$. We bound this term in the following way,
\[
II_{n}(N)\leq\frac{1}{n}\sum_{j=1}^{n}\sup_{m>N}\mathbb{E}_{0}[G_{m}^{2}%
(\xi_{j},\xi_{j-1})]\text{ .}%
\]
Notice that
\[
\sup_{m>N}\mathbb{E}_{0}[G_{m}^{2}(\xi_{k-1},\xi_{k})]\leq\mathbb{E}_{0}%
[\sup_{m>N}G_{m}^{2}(\xi_{k-1},\xi_{k})]
\]
Define $T(x):=\int_{S} \sup_{m>N}G_{m}^{2}(x,y)Q(x,dy)$. Then, by our
assumptions, $T\in L_{1}(\pi)$, and
\[
\frac{1}{n}\sum\nolimits_{k=1}^{n}\mathbb{E}_{0}[\sup_{m>N}G_{m}^{2}(\xi
_{k-1},\xi_{k})]=\frac{1}{n}\sum\nolimits_{k=1}^{n} Q^{k-1}T(\xi_{0}),
\]
and, by the ergodic theorem for $Q$, we obtain
\[
\lim_{n\rightarrow\infty}\frac{1}{n}\sum\nolimits_{k=1}^{n}\mathbb{E}_{0}%
[\sup_{m>N}G_{m}^{2}(\xi_{k-1},\xi_{k})]=\mathbb{E}[\sup_{m>N}G_{m}^{2}%
(\xi_{0},\xi_{1})]~\text{a.s.}
\]
Hence, for every $N\ge1$,
\[
\limsup_{n\rightarrow\infty}\frac{1}{n}\mathbb{E}_{0}[S_{n}-\mathbb{E}%
_{0}(S_{n})-M_{n}]^{2}\le\mathbb{E}[\sup_{m>N}G_{m}^{2}(\xi_{0},\xi_{1})]
~\text{a.s.}
\]
By (\ref{AS}) and dominated convergence theorem, $\mathbb{E}[\sup_{m>N}%
G_{m}^{2}(\xi_{0},\xi_{1})]\rightarrow0$ as $N\rightarrow\infty$ and
therefore
\[
\lim_{n\rightarrow\infty}\frac{1}{n}\mathbb{E}_{0}[S_{n}-\mathbb{E}_{0}%
(S_{n})-M_{n}]^{2}=0\text{ a.s. and in }\mathbb{L}_{1}\text{.}%
\]
Therefore (\ref{LimR}) is established.

To prove the quenched CLT we start from the approximation (\ref{LimR}) that we
combine with Theorem 3.1 in Billingsley (1999). It follows that for almost all
$x\in S$, the limiting distribution of $(S_{n}/\sqrt{n})_{n\geq1}$ is the same
as of $(M_{n}/\sqrt{n})_{n\geq1}$ under $\mathbb{P}^{x}$. Then, we shall use
the fact that an ergodic martingale with stationary differences satisfies the
CLT started at a point. For a complete and careful proof of this last fact we
direct to Derriennic and Lin (2001 a, page 520). $\lozenge$

\bigskip

We easily notice that:

\begin{remark}
\label{RlimAS}Condition (\ref{AS}) holds if and only if there is $H(\xi
_{0},\xi_{1})$ such that
\[
\mathbb{E}[\sup_{n>N}(H_{m}-H)^{2}(\xi_{0},\xi_{1})]\rightarrow0\text{ as
}N\rightarrow\infty.
\]

\end{remark}

Notice that actually we proved a stronger result.

\begin{remark}
Condition (\ref{AS}) in Proposition \ref{LimAS} can be replaced by
\[
H_{n}(\xi_{0},\xi_{1})\rightarrow H(\xi_{0},\xi_{1})\text{ in }\mathbb{L}_{2}%
\]
and
\[
\frac{1}{n}\sum_{j=1}^{n}\sup_{m>N}\mathbb{E}_{0}[(H_{m}-H)^{2}(\xi_{j}%
,\xi_{j-1})]\rightarrow0\text{ a.s.}%
\]

\end{remark}

Also, by the proof of Proposition \ref{LimAS} we can formulate a limit theorem
started from a point under centering.

\begin{remark}
Condition (\ref{AS}) in Proposition \ref{LimAS} implies and for almost all
$x\in S$%
\[
\frac{S_{n}-\mathbb{E}^{x}(S_{n})}{\sqrt{n}}\Longrightarrow|\sigma_{f}|N\text{
under }\mathbb{P}^{x}\text{ .}%
\]
where $\sigma_{f}^{2}=||H||_{2}^{2}.$
\end{remark}

For the sake of applications we give here two corollaries in terms of
conditional expectations of individual summands:

\begin{corollary}
Let $(\xi_{n})_{n\in\mathbb{Z}}$ be a stationary ergodic Markov chain and
$f\in${$\mathbb{L}_{0}^{2}(\pi)$} {satisfies condition (\ref{ASDL}) and in
addition }%
\begin{equation}
\sum_{k\geq1}||\mathbb{P}_{-k}(X_{0})||_{2}<\infty\text{ .} \label{Pr}%
\end{equation}
Then the conclusion of proposition \ref{LimAS} holds.
\end{corollary}

We know from Voln\'{y} (personal communication) that condition (\ref{Pr}) does
not imply condition {(\ref{ASDL}) and therefore this last condition cannot be
avoided.}

Next corollary provides a sufficient condition easy to verify. We adapt in its
proof an argument in Peligrad and Utev (2006, proof of Corollary 2).

\begin{corollary}
Let $(\xi_{n})_{n\in\mathbb{Z}}$ be a stationary ergodic Markov chain and
$f\in${$\mathbb{L}_{0}^{2}(\pi)$} {satisfies }
\begin{equation}
\sum_{k\geq1}\frac{||\mathbb{E}_{0}(X_{k})||_{2}}{\sqrt{k}}<\infty\text{ .}
\label{Mix}%
\end{equation}
Then the conclusion of proposition \ref{LimAS} holds.
\end{corollary}

\textbf{Proof}. By Cauchy-Schwarz, we have
\begin{gather*}
\sum_{k\geq1}||\mathbb{P}_{-k}(X_{0})||_{2} \le\sum_{n\ge0} 2^{n/2}%
(\sum_{k=2^{n}}^{2^{n+1}}||\mathbb{P}_{-k}(X_{0})||_{2}^{2})^{1/2} \le
\sum_{n\ge0} 2^{n/2} ||\mathbb{E}_{-{2^{n}}}(X_{0})||_{2},
\end{gather*}
and condition (\ref{Pr}) holds since $(||\mathbb{E}_{-{n}}(X_{0})||_{2})$ is
non increasing.
Then we notice that condition (\ref{Mix}) implies
\[
\sum_{k\geq1}\frac{|\mathbb{E}_{0}(X_{k})|}{\sqrt{k}}<\infty\text{ a.s. }%
\]
that further implies {(\ref{ASDL})} by Kronecker lemma.

\subsection{Almost sure results for reversible Markov chains}

In the context of normal Markov chains under condition (\ref{SN}), the
convergence in $\mathbb{L}_{2}$ below is known (Lemma 7.2 in Borodin and
Ibragimov, 1994). For reversible Markov chains we add here the almost sure convergence.

\begin{proposition}
\label{Key}Assume $(\xi_{n})_{n\in\mathbb{Z}}$ is a stationary ergodic
reversible Markov chain and $f\in${$\mathbb{L}_{0}^{2}(\pi)$ satisfies}
\eqref{SR2}. We have
\[
\mathbb{P}_{1}(S_{n})=H_{n}(\xi_{0},\xi_{1})\ \text{converges a.s. and in
}\mathbb{L}_{2}\text{ .}%
\]
Denote the limit by $H(\xi_{0},\xi_{1}).$ Moreover%
\[
\mathbb{E}_{0}[H(\xi_{0},\xi_{1})]=0\text{ \ a.s.}%
\]
and
\begin{equation}
\mathbb{E}(\sup_{k>n}[H_{k}-H]^{2}(\xi_{0},\xi_{1}))\rightarrow0\text{ as
}n\rightarrow\infty\text{ .} \label{lim0}%
\end{equation}

\end{proposition}

It is convenient to simplify the notation when no confusion is possible and we
further denote $H_{n}=H_{n}(\xi_{0},\xi_{1}).$ The proof will be divided into
two lemmas, allowing us to reduce the proof of the result for the sequences
$(H_{2n})$, $(H_{2^{n}})$ and then $(H_{2^{2^{n}}})$.

The proof follows a classical line, based on the dyadic chaining. We shall
estimate first, for every positive integers $\ m<n$, the expected values
$\mathbb{E}[H_{n}-H_{m}]^{2}$. Properties of conditional expectation and
Markov property give%
\[
\mathbb{E}[H_{n}-H_{m}]^{2}=\mathbb{E}[\sum_{m\leq i\leq n-1\text{ }}%
Q^{i}(f)(\xi_{1})]^{2}-\mathbb{E}[\sum_{m+1\leq i\leq n\text{ }}Q^{i}%
(f)(\xi_{0})]^{2}\text{ .}%
\]
Then, spectral calculus shows that (see page 166 in Borodin and Ibragimov,
1994)%
\begin{equation}
\mathbb{E}[H_{n}-H_{m}]^{2}=\int_{-1}^{1}(1-t^{2})[\sum_{k=n}^{m-1}t^{k}%
]^{2}\rho_{f}(dt)\text{.} \label{varH}%
\end{equation}

Using \eqref{varH} we will be able to obtain maximal inequalities along the
subsequences mentioned before.

\begin{lemma}
\label{L1}Assume $(\xi_{n})_{n\in\mathbb{Z}}$ is a stationary ergodic
reversible Markov chain and $f\in${$\mathbb{L}_{0}^{2}(\pi)$ satisfies}
\eqref{SR}. Then
\[
\sum_{d\geq0}\mathbb{E}\max_{2^{d}<n\leq2^{d+1}}[H_{2n}-H_{2^{d+1}}%
]^{2}<\infty.
\]

\end{lemma}

\textbf{ Proof.} We shall apply Lemma \ref{Wu} from Appendix. By (\ref{varH}),
for every $2^{d}\leq m<n\leq2^{d+1}-1,$ we estimate \
\begin{gather*}
\mathbb{E}[H_{2n}-H_{2m}]^{2}\leq\int_{-1}^{1}(1-t^{2})(\sum_{k=m}%
^{n-1}(t^{2k}+t^{2k+1}))^{2}\rho_{f}(dt)\\
\leq(n-m)^{2}\int_{-1}^{1}(1-t^{2})(1+t)^{2}t^{2^{d+1}}\rho_{f}(dt).
\end{gather*}
Take $T_{0}=H_{2^{d+1}},T_{1}=H_{2^{d+1}+2}\ldots,T_{2^{d}}=H_{2^{d+2}}$. For
every $0\leq r\leq d$, and $1\leq m\leq2^{d-r}$, we have
\[
\mathbb{E}[T_{2^{r}m}-T_{2^{r}(m-1)}]^{2}\leq2^{2r}\int_{-1}^{1}%
(1-t^{2})(1+t)^{2}t^{2^{d+1}}\rho_{f}(dt),
\]
which yields
\begin{gather*}
\mathbb{E}\max_{2^{d}\leq m<n\leq2^{d+1}-1}[H_{2n}-H_{2m}]^{2}\leq\int
_{-1}^{1}(1-t^{2})(1+t)^{2}t^{2^{d+1}}\rho_{f}(dt)(2^{d/2}\sum_{r=0}%
^{d}2^{r/2})^{2}\\
\leq2^{2d}\int_{-1}^{1}(1-t^{2})(1+t)^{2}t^{2^{d+1}}\rho_{f}(dt).
\end{gather*}
Now, for every $t\in(-1,1)$, we have
\[
\frac{1}{(1-t^{2})^{2}}=\sum_{k\geq1}(k-1)t^{2k}\geq\sum_{d\geq0}\sum
_{k=2^{d}}^{2^{d+1}-1}(k-1)t^{2k}\geq\sum_{d\geq0}2^{2(d-1)}t^{2^{d+1}}.
\]
Then, we obtain (notice that the cases $t\in\{-1,1\}$ do not cause any
problem)
\[
\sum_{d\geq0}\mathbb{E}\max_{2^{d}<n\leq2^{d+1}}[H_{2n}-H_{2^{d+1}}]^{2}%
\leq\int_{-1}^{1}\frac{1+t}{1-t}\rho_{f}(dt).
\]

$\lozenge$

\begin{lemma}
\label{L2}Assume $(\xi_{n})_{n\in\mathbb{Z}}$ is a stationary ergodic
reversible Markov chain and $f\in${$\mathbb{L}_{0}^{2}(\pi)$ satisfies}
\eqref{SR2}. Then $(H_{2^{n}})_{n\geq0}$ converges a.s. Moreover, we have
$\sup_{n\geq0}|H_{2^{n}}|\in L^{2}$.
\end{lemma}

\textbf{Proof.} We take $W_{n}=H_{2^{n}}$, and $\displaystyle g_{n}%
(t)=\sqrt{1-t^{2}}\sum_{k=2^{n}}^{2^{n+1}-1}t^{k}$ (which is positive on
$[-1,1]$), and apply Lemma \ref{lemme}. Hence we only need to prove that
condition \eqref{cond2} is satisfied, that is
\[
\int_{-1}^{1}(1-t^{2})\bigg(\sum_{n\geq1}(\log n)\sum_{k=2^{n}}^{2^{n+1}%
-1}t^{k}\bigg)^{2}\rho_{f}(dt)<\infty.
\]
Now,
\begin{gather*}
\sum_{n\geq1}(\log n)\sum_{k=2^{n}}^{2^{n+1}-1}t^{k}\leq C\sum_{n\geq1}%
\sum_{k=2^{n-1}}^{2^{n}}(\log\log(2k))(t^{2k}+t^{2k+1})\\
\leq C(1+t)\sum_{n\geq1}(\log\log(2k))t^{2k}\leq\tilde{C}\frac{(1+t)\log
^{+}|\log(1-t^{2})|}{1-t^{2}}\text{ },
\end{gather*}
where we used, for the last step a classical result, see e.g. Theorem 5, Ch.
XIV.5 of Feller (1971).

\smallskip Hence, it suffices that
\[
\int_{-1}^{1}\frac{(1+t)(\log^{+}|\log(1-t^{2})|)^{2}}{1-t}~\rho
_{f}(dt)<\infty\text{ ,}%
\]
which is clearly satisfied under our assumption. $\lozenge$

\bigskip

\textbf{Proof of Proposition \ref{Key}}. By combining the conclusions of Lemma
\ref{L1} and Lemma \ref{L2} we obtain
\[
\mathbb{E}(\sup_{k>n}[H_{2k}-H]^{2}(\xi_{0},\xi_{1}))\rightarrow0\text{ as
}n\rightarrow\infty\text{ .}%
\]
To finish the proof, it suffices to show that
\[
\sum_{n\ge1}|H_{n}-H_{n+1}|^{2}\in\mathbb{L}_{1}.
\]
Notice that $H_{n+1}-H_{n}=\mathbb{P}_{1}(X_{n+1}),$ where $\mathbb{P}%
_{1}(X_{n+1})$ is defined by (\ref{Defproj}). Since the projections are
orthogonal
\[
\sum_{j\geq0}\mathbb{E[P}_{1}(X_{j+1})]^{2}=\sum_{j\geq0}\mathbb{E[P}%
_{-j}(X_{0})]^{2}\leq\mathbb{E(}X_{0}^{2})\text{ .}%
\]
\smallskip

$\lozenge$

\bigskip

The following result is known from Derriennic and Lin (2001 b, Theorem 3.11,
ii). Their result is obtained under the condition $f\in Range(1-Q)^{1/2}$.
Kipnis and Varadhan (1986), in their discussion on page 4, proved that this
condition is equivalent to (\ref{var}).

\begin{lemma}
\label{DL}Assume $(\xi_{n})_{n\in\mathbb{Z}}$ is a stationary ergodic
reversible Markov chain and $f\in${$\mathbb{L}_{0}^{2}(\pi)$ satisfies}
\eqref{SR}. We have%
\[
\frac{\mathbb{E}_{0}(S_{n})}{\sqrt{n}}\rightarrow0\text{ almost surely and in
}\mathbb{L}_{2}\text{ .}%
\]

\end{lemma}

\textbf{Proof of Theorem \ref{REV}.}

\bigskip

To prove the validity of this theorem we shall apply Proposition \ref{LimAS}.
By Remark \ref{RlimAS} combined with Proposition \ref{Key} and Lemma \ref{DL}
we have the desired result.

\section{Application to random walks on compact groups}

In this section we shall apply our results to random walks on compact groups.

Let $\mathcal{X}$ be a compact abelian group, $\mathcal{A}$ a sigma algebra of
Borel subsets of $\mathcal{X}$ and $\pi$ the normalized Haar measure on
$\mathcal{X}$. The group operation is denoted by $+$. Let $\nu$ be a
probability measure on $(\mathcal{X},\mathcal{A)}$. The random walk on
$\mathcal{X}$ defined by $\nu$ is the Markov chain having the transition
function%
\[
(x,A)\rightarrow Q(x,A)=\nu(A-x)\text{ .}%
\]
The corresponding Markov operator denoted by $Q$ is defined by%
\[
(Qf)(x)=f\ast\nu(x)=\int_{\mathcal{X}}f(x+y)\nu(dy)\text{ .}%
\]
The Haar measure is invariant under $Q.$ We shall assume that $\nu$ is not
supported by a proper closed subgroup of $\mathcal{X},$ condition that is
equivalent to $Q$ being ergodic. In this context%
\[
(Q^{\ast}f)(x)=f\ast\nu^{\ast}(x)=\int_{\mathcal{X}}f(x-y)\nu(dy)\text{ ,}%
\]
where $\nu^{\ast}$ is the image of measure of $\nu$ by the map $x\rightarrow
-x.$ Thus $Q$ is symmetric on $\mathbb{L}_{2}(\pi)$ if and only if $\nu$ is
symmetric on $\mathcal{X}$, that is $\nu=\nu^{\ast}.$

The dual group of $\mathcal{X}$, denoted by $\mathcal{\hat{X}}$, is discrete.
Denote by $\hat{\nu}$ the Fourier transform of the measure $\nu,$ that is the
function%
\[
g\rightarrow\hat{\nu}(g)=\int_{\mathcal{X}}g(x)\nu(dx)\text{ }\ \text{with
}g\in\mathcal{\hat{X}}\text{ .}%
\]
A function $f\in${$\mathbb{L}^{2}(\pi)$} has the Fourier expansion%
\[
f=\sum\limits_{g\in\mathcal{\hat{X}}}\hat{f}(g)g\text{ .}%
\]
Ergodicity of $Q$ is equivalent to $\hat{\nu}(g)\neq1$ for any non-identity
$g\in\mathcal{\hat{X}}.$ By arguments in Borodin and Ibragimov (1994, ch. 4,
section 9) and also Derriennic and Lin (2001 a, Section 8) condition
(\ref{SR2}) takes the form%
\begin{equation}
\sum_{1\neq g\in\mathcal{\hat{X}}}\frac{|\hat{f}(g)|^{2}(\log^{+}|\log
|1-\hat{\nu}(g)|)^{2}}{|1-\hat{\nu}(g)|}<\infty\text{ .} \label{G1}%
\end{equation}

Combining these considerations with Theorem \ref{REV} we obtain the following result:

\begin{corollary}
Let $\nu$ be ergodic and symmetric on $\mathcal{X}$. If for $f$ in
{$\mathbb{L}_{0}^{2}(\pi)$ condition }(\ref{G1}) is satisfied then the central
limit theorem holds for the function $f$ and for the random walk generated by
$\nu,$ started at $x\in\mathcal{X}$ for $\pi-$almost every $x\in\mathcal{X}$.
\end{corollary}

A simple example is given by the following random walk on the one dimensional
torus which was considered in Paroux (1993) and Derriennic and Lin (2001 a and 2007).

On $\mathcal{X=}\mathbb{R}/\mathbb{Z}$ we take $\nu=\delta_{a}/2+\delta
_{-a}/2$ with $a\in\lbrack0,1]$ irrational. The random walk defined by $\nu$,
started at $x$, performs on the orbit of $x$ under the ergodic transformation
$T:x\rightarrow x+a$ $(mod$ $1)$, a simple random walk with the probability of
a step equal to $1/2$. The Fourier coefficients of $\nu$ are $\hat{\nu
}(n)=(e^{2i\pi na}+e^{-2i\pi na})/2$. For this example Derriennic and Lin
(2001) showed that $|1-\hat{\nu}(n)|\sim C\{na\}^{2},$ when the fractional
part of $na,$ denoted$\{na\},$ tends to $0$. With $\hat{f}(n)$ denoting the
$n^{th}$-Fourier coefficient of the function $f$ we have the following result:
If%
\[
\sum_{0\neq n\in Z}\frac{|\hat{f}(n)|^{2}(\log^{+}|\log\{na\}^{2}|)^{2}%
}{\{na\}^{2}}<\infty\text{ ,}%
\]
condition (\ref{G1}) is fulfilled and then, the CLT holds for the function $f$
and the random walk generated by $\nu$ started at $x$, for almost every
$x\in\mathbb{R}/\mathbb{Z}.$

If on $\mathcal{X=}\mathbb{R}/\mathbb{Z}$ we take $\nu=\delta_{a}%
/4+\delta_{-a}/4+\delta_{0}/2$, computations in Derriennic and Lin (2001 a)
show that condition (\ref{G1}) is satisfied as soon as
\[
\sum_{0\neq n\in Z}\frac{|\hat{f}(n)|^{2}(\log^{+}|\log(\sin^{2}(n\pi
\alpha)|)^{2}}{\sin^{2}(n\pi\alpha)}<\infty\text{ .}%
\]
Therefore, under this condition, the CLT holds for the function $f$ and the
random walk generated by $\nu,$ started at for almost every $x\in
\mathbb{R}/\mathbb{Z}.$

\section{Appendix}

The maximal inequality in this lemma is obtained by a chaining argument. It is
quoted from Wu (2007).

\begin{lemma}
\label{Wu}Let $p>1$ and $(T_{n})_{0\leq n\leq2^{d}}$ be random variables in
$\mathbb{L}_{2}$. Denote $T_{n}^{\ast}:=\max_{1\leq k\leq n}|T_{k}-T_{0}|$,
$1\leq n\leq2^{d}$. Then we have,
\begin{equation}
\Vert T_{2^{d}}^{\ast}\Vert_{2}\leq\sum_{r=0}^{d}[\sum_{m=1}^{2^{d-r}%
}\mathbb{E}(T_{2^{r}m}-T_{2^{r}(m-1)})^{2}]^{1/2}. \label{Wu1}%
\end{equation}

\end{lemma}

Next lemma provides a useful analytic tool to estimate maximum of partial sums
when the variance of partial sums is estimated as an integral. It is designed
to be applied directly in our proofs but the idea of proof should be useful in
a more general context. It is motivated by the paper M\'{o}ricz, Serfling and
Stout (1982).

\begin{lemma}
\label{lemme} Let $(W_{n})_{n\geq1}$ be a sequence of square-integrable
variables, and $(g_{n})$ a sequence of positive Borel functions on $[-1,1]$.
Assume that there exists a positive finite measure $\mu$ on $[-1,1]$, such
that
\begin{equation}
\mathbb{E}[W_{n}-W_{m}]^{2}\leq\int_{-1}^{1}[g_{m+1}(t)+\ldots+g_{n}%
(t)]^{2}\mu(dt)\quad\mbox{for every $1\le m < n$}. \label{cond}%
\end{equation}
If
\begin{equation}
\int_{-1}^{1}[\sum_{n\geq1}(\log n)g_{n}(t)]^{2}\mu(dt)<\infty\text{ },
\label{cond2}%
\end{equation}
then $(W_{n})_{n\geq1}$ converges almost surely and $\sup_{n\geq1}|W_{n}%
|\in\mathbb{L}_{2}$.
\end{lemma}

\textbf{Proof.}\newline Let $d$ be a positive integer. We shall apply Lemma
\ref{Wu} with $T_{0}=W_{2^{d}},T_{1}=W_{2^{d}+1},\ldots,T_{2^{d}}=W_{2^{d+1}}%
$. Notice that, by using \eqref{cond}, and the positivity of the functions
$g_{i}$ we have
\begin{gather*}
\sum_{m=1}^{2^{d-r}}\mathbb{E}[T_{2^{r}m}-T_{2^{r}(m-1)}]^{2}\leq\int_{-1}%
^{1}\sum_{m=1}^{2^{d-r}}[g_{2^{r}(m-1)+2^{d}+1}(t)\ldots+g_{2^{r}m+2^{d}%
}(t)]^{2}\mu(dt)\\
\leq\int_{-1}^{1}[\sum_{k=2^{d}+1}^{2^{d+1}}g_{k}(t)]^{2}\mu(dt)\text{ .}%
\end{gather*}
Hence, we obtain
\begin{equation}
\sum_{d\geq0}\mathbb{E}\max_{2^{d}+1\leq k\leq2^{d+1}}^{\ }[W_{k}-W_{2^{d}%
}]^{2}\leq\sum_{d\geq0}(d+1)^{2}\int_{-1}^{1}[\sum_{k=2^{d}+1}^{2^{d+1}}%
g_{k}(t)]^{2}\mu(dt).\label{RHS}%
\end{equation}
On the other hand, by Cauchy-Schwarz
\[
\lbrack\sum_{d\geq0}|W_{2^{d+1}}-W_{2^{d}}|]^{2}\leq\sum_{k\geq0}\frac
{1}{(k+1)^{2}}\sum_{d\geq0}(d+1)^{2}|W_{2^{d+1}}-W_{2^{d}}|^{2},
\]
and, by using \eqref{cond},
\begin{equation}
\mathbb{E}[\sum_{d\geq0}|W_{2^{d+1}}-W_{2^{d}}|]^{2}\leq\frac{\pi^{2}}{6}%
\sum_{d\geq0}(d+1)^{2}\int_{-1}^{1}[\sum_{k=2^{d}+1}^{2^{d+1}}g_{k}(t)]^{2}%
\mu(dt).\label{est}%
\end{equation}
By the positivity of the functions $g_{i}$, we see that condition
\eqref{cond2} implies the convergence of the right hand side of \eqref{RHS}
and (\ref{est}).

Hence $(W_{2^{n}})$ converges a.s. and $\max_{2^{d}+1\leq k\leq2^{d+1}}%
^{\ }|W_{k}-W_{2^{d}}|$ goes to 0 a.s., which yields the a.s. convergence of
$(W_{n})$. Moreover, for $2^{d}\leq n\leq2^{d+1}$,
\begin{gather*}
|W_{n}|\leq|W_{2^{d}}|+\max_{2^{d}+1\leq k\leq2^{d+1}}^{\ }|W_{k}-W_{2^{d}%
}|)\\
\leq|W_{1}|+\sum_{k=0}^{d-1}|W_{2^{k+1}}-W_{2^{k}}|+\max_{2^{d}+1\leq
k\leq2^{d+1}}^{\ }|W_{k}-W_{2^{d}}|,
\end{gather*}
and
\begin{gather*}
\sup_{n\geq1}W_{n}^{2}\leq3~\bigg(W_{1}^{2}+(\sum_{d\geq0}|W_{2^{k+1}%
}-W_{2^{k}}|)^{2}\\
+\sum_{d\geq0}\max_{2^{d}+1\leq k\leq2^{d+1}}^{\ }|W_{k}-W_{2^{d}}%
|^{2}\bigg)\in L_{1}.
\end{gather*}

$\lozenge$

\end{document}